\documentclass{amsart}

\usepackage{amssymb}
\usepackage{amsmath}
\usepackage{amsfonts}

\newtheorem{theorem}{Theorem}[section]
\theoremstyle{plain}

\newtheorem{corollary}[theorem]{Corollary}

\newtheorem{definition}[theorem]{Definition}

\newtheorem{lemma}[theorem]{Lemma}

\newtheorem{proposition}[theorem]{Proposition}
\newtheorem{remark}[theorem]{Remark}

\numberwithin{equation}{section}

\input{tcilatex}

\begin{document}
\title[Symmetry in $\mathcal{O}_{2}$]{Symmetry in the Cuntz Algebra on two
generators}
\author{Man-Duen Choi}
\address{Department of Mathematics\\
University of Toronto}
\email{choi@math.toronto.edu}
\author{Fr\'{e}d\'{e}ric Latr\'{e}moli\`{e}re}
\address{Department of Mathematics\\
University of Denver}
\email{frederic@math.du.edu}
\urladdr{http://www.math.du.edu/\symbol{126}frederic}
\date{February 22nd, 2010}

\begin{abstract}
We investigate the structure of the automorphism of $\mathcal{O}_{2}$ which
exchanges the two canonical isometries. Our main observation is that the
fixed point C*-subalgebra for this action is isomorphic to $\mathcal{O}_{2}$
and we detail the relationship between the crossed-product and fixed point
subalgebra.
\end{abstract}

\maketitle

\bigskip This paper studies the structure of the fixed point C*-algebra of
the action of $\mathbb{Z}_{2}$ which switches the canonical generators of
the Cuntz algebra $\mathcal{O}_{2}$. We show that both the
C*-crossed-product and the fixed point C*-algebra for this action are
*-isomorphic to $\mathcal{O}_{2}$.

This action is an example of an action of a finite group on a noncommutative
C*-algebra, and in general the structures associated to such actions can be
quite difficult to describe \cite{Rieffel80,Izumi04,Izumi04b}. To any action $\alpha$ of a finite group $G$ on a unital C*-algebra $A$, one can associate two new related C*-algebras:
the fixed point C*-algebra $A_{1}$ and the C*-crossed-product $A\rtimes
_{\alpha }G$ \cite{Zeller-Meier68}. By construction, $A_{1}$ is a
C*-subalgebra of $A$, while, if $G$ is Abelian, then $A$ is in fact the fixed point C*-subalgebra of $A\rtimes _{\alpha }G$ for the dual action of the Pontryagin dual of $G$ --- so $A$ is itself a subalgebra of $A\rtimes _{\alpha }G$. In \cite{Rosenberg79}
, Rosenberg shows that $A_{1}$ is *-isomorphic to a corner of $A\rtimes _{\alpha }G$, so that if $A\rtimes_\alpha G$ is simple, then it is Morita equivalent to $A_1$. In general, however, understanding the structure of $A_{1}$ or $
A\rtimes _{\alpha }G$ can be quite complex, as demonstrated for instance in 
\cite{Latremoliere06b}. In this paper, when $A$ is chosen to be $\mathcal{O}_{2}$ and the group is $\mathbb{Z}_{2}$, for the natural action swapping the
generators of $\mathcal{O}_{2}$, we obtain a complete picture of the
relative positions of these three C*-algebras, which we prove are all
*-isomorphic to $\mathcal{O}_{2}$.

We shall say that two isometries $S_{1}$ and $S_{2}$ on some Hilbert space
satisfy the Cuntz relation when:%
\begin{equation}
S_{1}S_{1}^{\ast }+S_{2}S_{2}^{\ast }=1\text{.}  \label{Cuntz}
\end{equation}

By \cite{Cuntz77}\cite[Theorem V.4.6 p. 147]{Davidson}, the Cuntz relation
defines, up to *-isomorphism, a unique simple C*-algebra denoted by $%
\mathcal{O}_{2}$. Moreover, by universality, there is a unique
*-automorphism $\sigma $ of $\mathcal{O}_{2}$ which satisfies:%
\begin{equation*}
\sigma (S_{1})=S_{2}\text{ and }\sigma (S_{2})=S_{1}\text{.}
\end{equation*}%
Since $\sigma ^{2}$ is the identity, we can define the C*-crossed-product $%
\mathcal{O}_{2}\rtimes _{\sigma }\mathbb{Z}_{2}$ as the universal C*-algebra
generated by two isometries $S_{1}$ and $S_{2}$ and a unitary $w$ such that $%
w^{2}=1$ and $wS_{1}=S_{2}w$ \cite{Zeller-Meier68}. We also can define the
fixed point C*-subalgebra $\left[ \mathcal{O}_{2}\right] _{1}$ of $\mathcal{O%
}_{2}$ as $\left\{ a\in \mathcal{O}_{2}:\sigma (a)=a\right\} $. It should also be noted that Izumi \cite[Example 5.7]{Izumi04} studied the action of $\mathbb{Z}_2$ on $\mathcal{O}_2$ given by $\sigma$ and proved that it has the Rohlin property. Thus, our examples fit in a larger family of ``classifiable actions'' in the sense of \cite{Izumi04}.

\bigskip In the first section of this paper, we show that $\left[ \mathcal{O}_{2}\right] _{1}$ is in fact *-isomorphic to $\mathcal{O}_{2}$. In the second section, we prove that $\mathcal{O}_{2}\rtimes _{\sigma }\mathbb{Z}_{2}$ is also *-isomorphic
to $\mathcal{O}_{2}$ and that $\sigma $ is not inner. We also establish that
in any representation of $\mathcal{O}_{2}$, the set of unitaries of order 2
exchanging the image of two generators is empty or a pair. In the third
section, we study the symmetry between the relations of $\mathcal{O}%
_{2}\rtimes _{\sigma }\mathbb{Z}_{2}$ with $\mathcal{O}_{2}$ on the one
hand, and $\mathcal{O}_{2}$ and $\left[ \mathcal{O}_{2}\right] _{1}$ on the
other hand. Section four deals with a description of the C*-crossed-product $%
\mathcal{O}_{2}\rtimes _{\sigma }\mathbb{Z}$, i.e. the universal C*-algebra
generated by a copy of $\mathcal{O}_{2}$ and a unitary $U$ such that $%
UaU^{\ast }=\sigma (a)$ for all $a\in \mathcal{O}_{2}$. We conclude this
paper with concrete representations of $\mathcal{O}_{2}$ on function spaces.

\section{Fixed Point C*-subalgebra}

\bigskip This section investigates the structure of the C*-algebra $\left[ 
\mathcal{O}_{2}\right] _{1}$ of fixed points of the automorphism $\sigma $.
We start with a simple preliminary result, which introduces a useful unitary
for our later purpose. We fix two isometries $S_{1}$ and $S_{2}$ satisfying
Relation (\ref{Cuntz}).

\begin{proposition}
\label{decomposition1}Let $S_{1}$ and $S_{2}$ be two isometries such that $%
S_{1}S_{1}^{\ast }+S_{2}S_{2}^{\ast }=1$ and $\sigma $ be the unique order
two automorphism of $\mathcal{O}_{2}=C^{\ast }\left( S_{1},S_{2}\right) $
such that $\sigma (S_{1})=S_{2}$. Let $U=S_{1}S_{1}^{\ast }-S_{2}S_{2}^{\ast
}$. Then $U$ is a unitary of $\mathcal{O}_{2}$ of order $2$. Let $\left[ 
\mathcal{O}_{2}\right] _{1}$ be the fixed point C*-subalgebra of $\mathcal{O}%
_{2}$ for $\sigma $. Then:%
\begin{equation*}
\mathcal{O}_{2}=\left[ \mathcal{O}_{2}\right] _{1}\oplus \left[ \mathcal{O}%
_{2}\right] _{1}U
\end{equation*}%
with $U$ a unitary of order $2$. Moreover, with this decomposition, if $%
a=a_{1}+a_{2}U$ then $\sigma (a)=a_{1}-a_{2}U$. Note that $\oplus $ is the
direct sum for Banach spaces, not between algebras, since $\left[ \mathcal{O}%
_{2}\right] _{1}U$ is not an algebra for the multiplication of $\mathcal{O}%
_{2}$.
\end{proposition}

\begin{proof}
For all $a\in \mathcal{O}_{2}$ we have $a=a_{1}+a_{-1}$ with $a_{\varepsilon
}=\frac{1}{2}\left( a+\varepsilon \sigma (a)\right) ,$ so that $\sigma
(a_{\varepsilon })=\varepsilon a_{\varepsilon }$ for $\varepsilon \in
\left\{ -1,1\right\} $. Let $\left[ \mathcal{O}_{2}\right] _{-1}$ be the
space of elements $a\in \mathcal{O}_{2}$ such that $\sigma (a)=-a$. It is
then immediate that $\mathcal{O}_{2}=\left[ \mathcal{O}_{2}\right]
_{1}\oplus \left[ \mathcal{O}_{2}\right] _{-1}$. Now, by construction, $%
U=U^{\ast }$ and:%
\begin{equation*}
U^{2}=\left( S_{1}S_{1}^{\ast }-S_{2}S_{2}^{\ast }\right) \left(
S_{1}S_{1}^{\ast }-S_{2}S_{2}^{\ast }\right) =S_{1}S_{1}^{\ast
}+S_{2}S_{2}^{\ast }=1\text{,}
\end{equation*}%
so $U$ is an order $2$ unitary. Moreover $\sigma (U)=-U$. Thus $a\in \left[ 
\mathcal{O}_{2}\right] _{-1}$ if and only if $aU\in \left[ \mathcal{O}_{2}%
\right] _{1}$. Hence our decomposition is proven. An immediate computation
shows that $\sigma $ is indeed implemented as shown.
\end{proof}

\bigskip We now start the process to identify $\left[ \mathcal{O}_{2}\right]
_{1}$. Our proof will exhibit a specific and interesting choice of
generators for $\left[ \mathcal{O}_{2}\right] _{1}$, and for clarity of
exposition it will be useful to keep track of the generators of the many
*-isomorphic copies of $\mathcal{O}_{2}$ we will encounter in our proof. We
start with the following lemmas:

\begin{lemma}
\label{FixedPointLemma1}Let $S_{1}$ and $S_{2}$ be two isometries such that $%
S_{1}S_{1}^{\ast }+S_{2}S_{2}^{\ast }=1$, so that $\mathcal{O}_{2}=C^{\ast
}\left( S_{1},S_{2}\right) $. We define the following elements in $%
M_{2}\left( C^{\ast }\left( S_{1},S_{2}\right) \right) $:%
\begin{equation*}
T_{1}=\frac{1}{\sqrt[2]{2}}\left[ 
\begin{array}{cc}
S_{1} & S_{2} \\ 
S_{1} & S_{2}%
\end{array}%
\right] \text{ and }T_{2}=\frac{1}{\sqrt[2]{2}}\left[ 
\begin{array}{cc}
S_{1} & -S_{2} \\ 
-S_{1} & S_{2}%
\end{array}%
\right] \text{.}
\end{equation*}%
Then $T_{1}^{\ast }T_{1}=T_{2}^{\ast }T_{2}=T_{1}T_{1}^{\ast
}+T_{2}T_{2}^{\ast }=1$ in $M_{2}\left( C^{\ast }\left( S_{1},S_{2}\right)
\right) $. Thus, by universality and simplicity of $\mathcal{O}_{2}$, the
C*-algebra $C^{\ast }\left( T_{1},T_{2}\right) $ is *-isomorphic to $%
\mathcal{O}_{2}$. On the other hand:%
\begin{equation*}
C^{\ast }\left( T_{1},T_{2}\right) =\left\{ \left[ 
\begin{array}{cc}
A_{1} & A_{2} \\ 
\sigma (A_{2}) & \sigma (A_{1})%
\end{array}%
\right] :A_{1},A_{2}\in C^{\ast }\left( S_{1},S_{2}\right) \right\}
\end{equation*}%
where $\sigma $ is the unique order-2 automorphism of $C^{\ast }\left(
S_{1},S_{2}\right) $ such that $\sigma (S_{1})=S_{2}$.
\end{lemma}

\begin{proof}
Note that $\sigma (S_{2})=\sigma \left( \sigma \left( S_{1}\right) \right)
=S_{1}$ by assumption on $\sigma $. Now, we observe that:%
\begin{equation*}
T_{1}=\frac{1}{\sqrt[2]{2}}\left( \left[ 
\begin{array}{cc}
S_{1} & 0 \\ 
0 & \sigma (S_{1})%
\end{array}%
\right] +\left[ 
\begin{array}{cc}
0 & S_{2} \\ 
\sigma (S_{2}) & 0%
\end{array}%
\right] \right)
\end{equation*}%
and%
\begin{equation*}
T_{2}=\frac{1}{\sqrt[2]{2}}\left( \left[ 
\begin{array}{cc}
S_{1} & 0 \\ 
0 & \sigma (S_{1})%
\end{array}%
\right] -\left[ 
\begin{array}{cc}
0 & S_{2} \\ 
\sigma (S_{2}) & 0%
\end{array}%
\right] \right)
\end{equation*}%
so $C^{\ast }\left( T_{1},T_{2}\right) =C^{\ast }\left( \left[ 
\begin{array}{cc}
S_{1} & 0 \\ 
0 & \sigma (S_{1})%
\end{array}%
\right] ,\left[ 
\begin{array}{cc}
0 & S_{2} \\ 
\sigma (S_{2}) & 0%
\end{array}%
\right] \right) $.

On the other hand, we also have:%
\begin{equation*}
\left[ 
\begin{array}{cc}
0 & 1 \\ 
1 & 0%
\end{array}%
\right] =\left[ 
\begin{array}{cc}
S_{1} & 0 \\ 
0 & \sigma (S_{1})%
\end{array}%
\right] \left[ 
\begin{array}{cc}
0 & S_{2} \\ 
\sigma (S_{2}) & 0%
\end{array}%
\right] ^{\ast }+\left[ 
\begin{array}{cc}
0 & S_{2} \\ 
\sigma (S_{2}) & 0%
\end{array}%
\right] \left[ 
\begin{array}{cc}
S_{1} & 0 \\ 
0 & \sigma (S_{1})%
\end{array}%
\right] ^{\ast }
\end{equation*}%
so $\Xi =\left[ 
\begin{array}{cc}
0 & 1 \\ 
1 & 0%
\end{array}%
\right] \in C^{\ast }\left( T_{1},T_{2}\right) $ and in fact:%
\begin{equation*}
\left[ 
\begin{array}{cc}
S_{2} & 0 \\ 
0 & \sigma (S_{2})%
\end{array}%
\right] =\left[ 
\begin{array}{cc}
0 & S_{2} \\ 
\sigma (S_{2}) & 0%
\end{array}%
\right] \Xi \text{.}
\end{equation*}%
Thus we conclude:%
\begin{equation}
C^{\ast }\left( T_{1},T_{2}\right) =C^{\ast }\left( \left[ 
\begin{array}{cc}
S_{1} & 0 \\ 
0 & \sigma (S_{1})%
\end{array}%
\right] ,\left[ 
\begin{array}{cc}
S_{2} & 0 \\ 
0 & \sigma (S_{2})%
\end{array}%
\right] ,\left[ 
\begin{array}{cc}
0 & 1 \\ 
1 & 0%
\end{array}%
\right] \right) \text{.}  \label{FixedPointLemma1Eq1}
\end{equation}%
We note that:%
\begin{equation*}
C^{\ast }\left( \left[ 
\begin{array}{cc}
S_{1} & 0 \\ 
0 & \sigma (S_{1})%
\end{array}%
\right] ,\left[ 
\begin{array}{cc}
S_{2} & 0 \\ 
0 & \sigma (S_{2})%
\end{array}%
\right] \right) =\left\{ \left[ 
\begin{array}{cc}
A & 0 \\ 
0 & \sigma (A)%
\end{array}%
\right] :A\in C^{\ast }\left( S_{1},S_{2}\right) \right\} \text{.}
\end{equation*}%
Thus, if $A_{1},A_{2}\in C^{\ast }\left( S_{1},S_{2}\right) $ then:%
\begin{equation*}
\left[ 
\begin{array}{cc}
A_{1} & A_{2} \\ 
\sigma (A_{2}) & \sigma (A_{1})%
\end{array}%
\right] =\left[ 
\begin{array}{cc}
A_{1} & 0 \\ 
0 & \sigma (A_{1})%
\end{array}%
\right] +\left[ 
\begin{array}{cc}
A_{2} & 0 \\ 
0 & \sigma (A_{2})%
\end{array}%
\right] \Xi
\end{equation*}%
is in $C^{\ast }\left( T_{1},T_{2}\right) $. Conversely, if we write $%
D_{2}\left( C^{\ast }\left( S_{1},S_{2}\right) \right) $ the algebra of
diagonal matrices in $M_{2}\left( C^{\ast }\left( S_{1},S_{2}\right) \right) 
$ then:%
\begin{equation*}
M_{2}\left( C^{\ast }\left( S_{1},S_{2}\right) \right) =D_{2}\left( C^{\ast
}\left( S_{1},S_{2}\right) \right) \oplus D_{2}\left( C^{\ast }\left(
S_{1},S_{2}\right) \right) \Xi \text{.}
\end{equation*}
Hence any element of $C^{\ast }\left( T_{1},T_{2}\right) $ must be of the
desired form from Equation (\ref{FixedPointLemma1Eq1}), which concludes our
lemma.
\end{proof}

\begin{lemma}
\label{FixedPointLemma2}We use the notations of Lemma (\ref{FixedPointLemma1}%
). Let $Z=\left[ 
\begin{array}{cc}
1 & 0 \\ 
0 & -1%
\end{array}%
\right] \in M_{2}\left( C^{\ast }\left( S_{1},S_{2}\right) \right) $ and let:%
\begin{equation*}
\tau :X\in C^{\ast }\left( T_{1},T_{2}\right) \mapsto ZXZ\text{.}
\end{equation*}%
Then $\tau $ is an order-2 automorphism on $C^{\ast }\left(
T_{1},T_{2}\right) $ such that $\tau (T_{1})=T_{2}$ and the fixed point
C*-algebra of $\tau $ is given by:%
\begin{equation*}
\left\{ \left[ 
\begin{array}{cc}
A & 0 \\ 
0 & \sigma (A)%
\end{array}%
\right] :A\in C^{\ast }\left( S_{1},S_{2}\right) \right\} \text{.}
\end{equation*}
\end{lemma}

\begin{proof}
The fixed point algebra of $C^{\ast }\left( T_{1},T_{2}\right) $ for $\tau $
is given by:%
\begin{equation*}
\left\{ a+\tau (a):a\in C^{\ast }\left( T_{1},T_{2}\right) \right\} \text{,}
\end{equation*}
so this lemma follows from an immediate computation and Lemma (\ref%
{FixedPointLemma1}).
\end{proof}

\begin{theorem}
\label{FixedPointAlgebraThm}Let $S_{1}$ and $S_{2}$ be two isometries such
that $S_{1}S_{1}^{\ast }+S_{2}S_{2}^{\ast }=1$ and $\sigma $ be the unique
order two automorphism of $\mathcal{O}_{2}=C^{\ast }\left(
S_{1},S_{2}\right) $ such that $\sigma (S_{1})=S_{2}$. Let:%
\begin{align*}
T& =\frac{1}{\sqrt[2]{2}}\left( S_{1}+S_{2}\right) \text{, } \\
U& =S_{1}S_{1}^{\ast }-S_{2}S_{2}^{\ast }\text{ and } \\
V& =UTU=\frac{1}{\sqrt[2]{2}}\left( S_{1}-S_{2}\right) \left(
S_{1}S_{1}^{\ast }-S_{2}S_{2}^{\ast }\right) \text{.}
\end{align*}%
Then the fixed point C*-algebra $\left[ \mathcal{O}_{2}\right] _{1}$\ for $%
\sigma $ is $C^{\ast }\left( T,V\right) $ and is *-isomorphic to $\mathcal{O}%
_{2}$.
\end{theorem}

\begin{proof}
We shall use the notations of Lemma (\ref{FixedPointLemma1}). First, let $%
\Phi :C^{\ast }\left( S_{1},S_{2}\right) \longrightarrow C^{\ast }\left(
T_{1},T_{2}\right) $ be the unique *-epimorphism defined by universality
with $\Phi (S_{j})=T_{j}$ ($j=1,2$). Since $C^{\ast }\left(
S_{1},S_{2}\right) $ is simple, $\Phi $ is a *-isomorphism. Moreover, by
construction, $\Phi \circ \sigma =\tau \circ \Phi $. Therefore, the fixed
point C*-algebra for $\sigma $ is *-isomorphic to the fixed point C*-algebra
for $\tau $.

Now, the fixed point C*-algebra for $\tau $ is given by Lemma (\ref%
{FixedPointLemma2}) as:%
\begin{equation*}
\left\{ \left[ 
\begin{array}{cc}
A & 0 \\ 
0 & \sigma (A)%
\end{array}%
\right] :A\in C^{\ast }\left( S_{1},S_{2}\right) \right\}
\end{equation*}%
so it is the C*-algebra generated by $R_{1}=\left[ 
\begin{array}{cc}
S_{1} & 0 \\ 
0 & S_{2}%
\end{array}%
\right] $ and $R_{2}=\left[ 
\begin{array}{cc}
S_{2} & 0 \\ 
0 & S_{1}%
\end{array}%
\right] $, which are two isometries satisfying the Cuntz relation. So the
fixed point C*-algebra for $\tau $ (hence for $\sigma $) is *-isomorphic to $%
\mathcal{O}_{2}$.

On the other hand, we have the relation:%
\begin{equation}
R_{1}=\frac{1}{\sqrt[2]{2}}\left( T_{1}+T_{2}\right)  \label{FixedPointEq1}
\end{equation}%
Moreover, if $Y=\Phi (U)$ then:%
\begin{eqnarray*}
Y &=&T_{1}T_{1}^{\ast }-T_{2}T_{2}^{\ast } \\
&=&\frac{1}{2}\left( \left[ 
\begin{array}{cc}
S_{1} & S_{2} \\ 
S_{1} & S_{2}%
\end{array}%
\right] \left[ 
\begin{array}{cc}
S_{1}^{\ast } & S_{1}^{\ast } \\ 
S_{2}^{\ast } & S_{2}^{\ast }%
\end{array}%
\right] -\left[ 
\begin{array}{cc}
S_{1} & -S_{2} \\ 
-S_{1} & S_{2}%
\end{array}%
\right] \left[ 
\begin{array}{cc}
S_{1}^{\ast } & -S_{1}^{\ast } \\ 
-S_{2}^{\ast } & S_{2}^{\ast }%
\end{array}%
\right] \right) \\
&=&\left[ 
\begin{array}{cc}
0 & 1 \\ 
1 & 0%
\end{array}%
\right]
\end{eqnarray*}%
so we obtain the relation:%
\begin{equation}
R_{2}=YR_{1}Y\text{.}  \label{FixedPointEq2}
\end{equation}

This concludes our proof after application of $\Phi ^{-1}$ to Relations (\ref%
{FixedPointEq1}) and (\ref{FixedPointEq2}).
\end{proof}

\bigskip Thus, Proposition (\ref{decomposition1}) can now be restated in the
following manner:\ $\mathcal{O}_{2}\mathcal{\ }$is *-isomorphic to $\mathcal{%
O}_{2}\oplus \mathcal{O}_{2}U$, where $\sigma (a\oplus bU)=a-bU$ for any $%
a,b\in \mathcal{O}_{2}$. Moreover, we have a pair of natural generators for $%
\left[ \mathcal{O}_{2}\right] _{1}$. It is natural to ask whether this
decomposition, in fact, is a mean to recognize $\mathcal{O}_{2}$ as a
crossed-product of an action on $\mathcal{O}_{2}$ implemented by $\limfunc{Ad%
}U$, and $\sigma $ can then be seen as the dual action of $\mathbb{Z}_{2}$
on this crossed-product. We note that $\limfunc{Ad}U$ does swap the
generators $T$ and $V$ of $\left[ \mathcal{O}_{2}\right] _{1}$ with the
notation of Theorem (\ref{FixedPointAlgebraThm}). The next two sections will
make precise these informal observations. We start with a study of the
structure of the C*-crossed-product $\mathcal{O}_{2}\rtimes _{\sigma }%
\mathbb{Z}_{2}$.

\section{C*-Crossed-Product}

We first observe that the C*-crossed-product $\mathcal{O}_{2}\rtimes
_{\sigma }\mathbb{Z}_{2}$ is in fact *-isomorphic to $\mathcal{O}_{2}$:

\begin{theorem}
\label{CrossedProductThm}Let $S_{1}$ and $S_{2}$ be two isometries such that 
$S_{1}S_{1}^{\ast }+S_{2}S_{2}^{\ast }=1$ and $\sigma $ be the unique order
two automorphism of $\mathcal{O}_{2}=C^{\ast }\left( S_{1},S_{2}\right) $
such that $\sigma (S_{1})=S_{2}$. Then:%
\begin{equation*}
\mathcal{O}_{2}\rtimes _{\sigma }\mathbb{Z}_{2}\text{ is *-isomorphic to }%
\mathcal{O}_{2}\text{.}
\end{equation*}
\end{theorem}

\begin{proof}
Let $W$ be the canonical unitary in $\mathcal{O}_{2}\rtimes _{\sigma }%
\mathbb{Z}_{2}$ such that $WS_{1}W=S_{2}$ and $W^{2}=1$. Then:%
\begin{eqnarray*}
S_{1}S_{1}^{\ast }W+\left( S_{1}S_{1}^{\ast }W\right) ^{\ast }
&=&S_{1}S_{1}^{\ast }W+WS_{1}S_{1}^{\ast } \\
&=&S_{1}S_{1}^{\ast }W+WS_{1}S_{1}^{\ast }W^{2} \\
&=&S_{1}S_{1}^{\ast }W+S_{2}S_{2}^{\ast }W \\
&=&W\text{.}
\end{eqnarray*}%
Hence, $W\in C^{\ast }\left( S_{1},WS_{1}\right) \subseteq \mathcal{O}%
_{2}\rtimes _{\sigma }\mathbb{Z}_{2}$. Since $\mathcal{O}_{2}\rtimes
_{\sigma }\mathbb{Z}_{2}$ is generated by $S_{1},S_{2}$ and $W$ and $%
S_{2}=WS_{1}W\in C^{\ast }\left( S_{1},WS_{1}\right) $ so $\mathcal{O}%
_{2}\rtimes _{\sigma }\mathbb{Z}_{2}=C^{\ast }\left( S_{1},WS_{1}\right) $.

On the other hand:%
\begin{equation*}
\left( WS_{1}\right) ^{\ast }WS_{1}=S_{1}^{\ast }W^{2}S_{1}=S_{1}^{\ast
}S_{1}=1
\end{equation*}%
and%
\begin{equation*}
S_{1}S_{1}^{\ast }+WS_{1}\left( WS_{1}\right) ^{\ast }=S_{1}S_{1}^{\ast
}+WS_{1}S_{1}^{\ast }W=S_{1}S_{1}^{\ast }+S_{2}S_{2}^{\ast }=1\text{.}
\end{equation*}%
Therefore, $C^{\ast }\left( S_{1},WS_{1}\right) $ is *-isomorphic to $%
\mathcal{O}_{2}$.
\end{proof}

\bigskip We can provide more details on the structure of the automorphism $%
\sigma $.

\begin{proposition}
\label{notinner}Let $S_{1}$ and $S_{2}$ be two isometries such that $%
S_{1}S_{1}^{\ast }+S_{2}S_{2}^{\ast }=1$ and $\sigma $ be the unique order
two automorphism of $\mathcal{O}_{2}=C^{\ast }\left( S_{1},S_{2}\right) $
such that $\sigma (S_{1})=S_{2}$. Then $\sigma $ is not inner.
\end{proposition}

\begin{proof}
Let $\mathcal{H}=l^{2}\left( \mathbb{N}\right) $ whose canonical Hilbert
basis is denoted by $\left( e_{n}\right) _{n\in \mathbb{N}}$ (namely, $%
\left( e_{n}\right) _{m}$ is $0$ unless $n=m$, when it is $1$). We define:%
\begin{equation*}
T_{1}e_{n}=e_{2n}\text{ and }T_{2}e_{n}=e_{2n+1}\text{.}
\end{equation*}%
Then note that $T_{2}$ has no eigenvector while $T_{1}e_{0}=e_{0}$. Hence, $%
T_{1}$ and $T_{2}$ are not unitarily equivalent in $\mathcal{H}$. Yet, it is
immediate that $T_{1}$ and $T_{2}$ are isometries which satisfy $%
T_{1}T_{1}^{\ast }+T_{2}T_{2}^{\ast }=1$. Therefore, there exists a (unique)
*-homomorphism $\varphi $ from $\mathcal{O}_{2}$ onto $C^{\ast }\left(
T_{1},T_{2}\right) $ with $\varphi (S_{j})=T_{j}$ for $j=1,2$, and since $%
\mathcal{O}_{2}$ is simple, $\varphi $ is in fact a *-monomorphism. Now, if $%
\sigma $ was inner, then there would exists some unitary $u\in \mathcal{O}%
_{2}$ such that $uS_{1}u^{\ast }=\sigma (S_{1})=S_{2}$. This would imply
that $\varphi (u)T_{1}\varphi (u)^{\ast }=T_{2}$ with $\varphi (u)$ a
unitary. This is a contradiction.
\end{proof}

\bigskip We can use Proposition (\ref{notinner}) to see that, if we can find
a covariant representation of $\mathcal{O}_{2}$, then the representation of $%
\mathbb{Z}_{2}$ is unique up to a sign.

\begin{proposition}
\label{uniqueness}Let $S_{1}$ and $S_{2}$ be two isometries such that $%
S_{1}S_{1}^{\ast }+S_{2}S_{2}^{\ast }=1$. Let $u$ and $w$ be two unitaries
such that $C^{\ast }\left( S_{1},S_{2},u\right) \subseteq C^{\ast }\left(
S_{1},S_{2},w\right) $, and such that $u^{2}=w^{2}=1$ with $uS_{1}=S_{2}u$
and $wS_{1}=S_{2}w$. Then $u=w$ or $u=-w$.
\end{proposition}

\begin{proof}
By assumption, $C^{\ast }\left( S_{1},S_{2}\right) $ is *-isomorphic to $%
\mathcal{O}_{2}$ since $\mathcal{O}_{2}$ is simple and universal for the
given property. Moreover, by universality, there exists a (unique)
*-morphism $\varphi :\mathcal{O}_{2}\rtimes _{\sigma }\mathbb{Z}%
_{2}\twoheadrightarrow C^{\ast }\left( S_{1},S_{2},w\right) $. Since $%
\mathcal{O}_{2}\rtimes _{\sigma }\mathbb{Z}_{2}$ is $\mathcal{O}_{2}$, hence
simple, $\varphi $ is an isomorphism.

We can use $\varphi $ to show that $C^{\ast }\left( S_{1},S_{2},w\right) =%
\mathcal{O}_{2}\oplus \mathcal{O}_{2}w$. Let us now write $u=a+bw$ with $%
a,b\in \mathcal{O}_{2}$. Then, for $j=1,2$ and by assumption, $%
uwS_{j}=S_{j}uw$, so:%
\begin{equation*}
\left( b+aw\right) S_{j}=S_{j}\left( b+aw\right) \text{.}
\end{equation*}%
Since $\mathcal{O}_{2}$ and $\mathcal{O}_{2}w$ are complementary spaces, we
conclude that $bS_{j}=S_{j}b$ and $awS_{j}=S_{j}aw$ (note that $awS_{j}=aS_{j^{\prime }}w$ with $aS_{j^{\prime }}\in \mathcal{O}_{2}$). Thus 
$b$ is in the center of $\mathcal{O}_{2}$ and thus is scalar. On the other
hand, we have:%
\begin{equation*}
\left\{ 
\begin{array}{c}
aS_{2}w=S_{1}aw \\ 
aS_{1}w=S_{2}aw%
\end{array}%
\right. \text{ so }\left\{ 
\begin{array}{c}
aS_{2}=S_{1}a \\ 
aS_{1}=S_{2}a%
\end{array}%
\right. \text{.}
\end{equation*}%
Consequently, $a^{2}$ commutes with $S_{1}$ and $S_{2}$ so it is central in $%
\mathcal{O}_{2}$, hence again $a^{2}$ is scalar, say $\lambda \in \mathbb{C}$%
. Now, since $u$ is normal and $w$ is normal, so are $a$ and $b$ (again,
since $\mathcal{O}_{2}$ and $\mathcal{O}_{2}w$ are complementary spaces).
Assume $a\not=0$. Then $v=\mu a$, where $\mu ^{2}=\lambda ^{-1}$, is a
unitary of order $2$ in $C^{\ast }\left( S_{1},S_{2}\right) $\ which
satisfies $vS_{j}=S_{j^{\prime }}v$. By Proposition (\ref{notinner}), this
is not possible. Hence, $a^{2}=0$ and so $a=0$ as $a$ normal. Thus $u=bw$
with $b$ scalar, and since $1=u^{2}=w^{2}$ we conclude that $b\in \left\{
-1,1\right\} $.
\end{proof}

\begin{remark}
We can recover the well-known fact that $\mathcal{O}_{2}=M_{2}\left( 
\mathcal{O}_{2}\right) $ as proven in \cite{Choi79}. Indeed, let us use the
notations of Lemma (\ref{FixedPointLemma1}). Then:%
\begin{equation*}
\mathcal{O}_{2}=C^{\ast }\left( T_{1},T_{2}\right) =\left\{ \left[ 
\begin{array}{cc}
a & b \\ 
\sigma (b) & \sigma (a)%
\end{array}%
\right] :a,b\in \mathcal{O}_{2}\right\}
\end{equation*}%
and, if $Z\in M_{2}\left( \mathcal{O}_{2}\right) $ with $Z=\left[ 
\begin{array}{cc}
1 & 0 \\ 
0 & -1%
\end{array}%
\right] $ then:%
\begin{equation*}
M_{2}\left( \mathcal{O}_{2}\right) =C^{\ast }\left( T_{1},T_{2},Z\right) 
\text{.}
\end{equation*}%
On the other hand, $C^{\ast }\left( T_{1},T_{2},Z\right) =C^{\ast }\left(
T_{1},T_{2}\right) \rtimes _{\eta }\mathbb{Z}_{2}$ where $\eta (T_{1})=T_{2}$
is of order two. Indeed, by universality, $C^{\ast }\left(
T_{1},T_{2},Z\right) $ is a quotient of $C^{\ast }\left( T_{1},T_{2}\right)
\rtimes _{\eta }\mathbb{Z}_{2}$, yet the later is $\mathcal{O}_{2}$ by
Theorem (\ref{CrossedProductThm}) so it is simple. Moreover, Theorem (\ref%
{CrossedProductThm}) provides us with a natural pair of generators for $%
M_{2}\left( \mathcal{O}_{2}\right) $.
\end{remark}

\bigskip Now, we wish to see that in some way, the embedding of $\mathcal{O}%
_{2}$ as the fixed point algebra for $\sigma $ in $\mathcal{O}_{2}$ or the
embedding of $\mathcal{O}_{2}$ into $\mathcal{O}_{2}\rtimes _{\sigma }%
\mathbb{Z}_{2}$ are the same. The following section formalizes this
statement.

\section{A doubly infinite Sequence of self-similar $\mathcal{O}_{2}$
embeddings}

\bigskip The C*-algebra $\mathcal{O}_{2}$ embeds into itself as a fixed
point sub-C*-algebra for $\sigma $ or as a subalgebra of its
crossed-product. The second embedding can be seen as embedding a fixed point
for the dual action to $\sigma $. In our case, these two embeddings are the
same, as shown in the following proposition.

\begin{proposition}
\label{Tau}There exists a *-isomorphism $\tau :\mathcal{O}_{2}\rtimes
_{\sigma }\mathbb{Z}_{2}\longrightarrow \mathcal{O}_{2}$ such that $\tau (%
\mathcal{O}_{2})$ is the fixed point C*-algebra of $\sigma $.
\end{proposition}

\begin{proof}
In Theorem\ (\ref{CrossedProductThm}), we showed that $\mathcal{O}%
_{2}\rtimes _{\sigma }\mathbb{Z}_{2}=C^{\ast }\left( S_{1},S_{1}W\right) $.
It follows that $\mathcal{O}_{2}\rtimes _{\sigma }\mathbb{Z}_{2}=C^{\ast
}\left( \frac{S_{1}+S_{1}W}{\sqrt[2]{2}},\frac{S_{1}-S_{1}W}{\sqrt[2]{2}}%
\right) $ and a direct computation shows that $B_{1}=\frac{S_{1}+S_{1}W}{%
\sqrt[2]{2}}$ and $B_{2}=\frac{S_{1}-S_{1}W}{\sqrt[2]{2}}$ are isometries
satisfying $B_{1}B_{1}^{\ast }+B_{2}B_{2}^{\ast }=1$. By universality of $%
\mathcal{O}_{2}$ there exists a unique *-monomorphism $\tau :\mathcal{O}%
_{2}\rtimes _{\sigma }\mathbb{Z}_{2}\longrightarrow \mathcal{O}_{2}$ defined
by $\tau (B_{i})=S_{i}$ for $i=1,2$. Now, $\tau (S_{1})=\tau \left( \frac{%
\sqrt[2]{2}}{2}\left( B_{1}+B_{2}\right) \right) =S_{1}+S_{2}=T$ and $\tau
(S_{2})=\tau \left( WS_{1}W\right) =UTU$ where $U=S_{1}S_{1}^{\ast
}-S_{2}S_{2}^{\ast }$ and $T=S_{1}+S_{2}$ following the notations of Theorem
(\ref{FixedPointAlgebraThm}). Hence $\tau $ maps $\mathcal{O}_{2}\subseteq 
\mathcal{O}_{2}\rtimes _{\sigma }\mathbb{Z}_{2}$ onto the fixed point
subalgebra of $\mathcal{O}_{2}$ for $\sigma $.
\end{proof}

\bigskip Now, we can give a somewhat more detailed picture of the embeddings
by seeing how we can construct a double infinite sequence of identical
embeddings of $\mathcal{O}_{2}$ into itself using the crossed-product
construction. More precisely, suppose we are given a copy of $\mathcal{O}%
_{2} $ generated by two isometries $r_{n}$ and $t_{n}$ with $%
r_{n}r_{n}^{\ast }+t_{n}t_{n}^{\ast }=1$ where $n\in \mathbb{Z}$ arbitrary.
Then we can define $\sigma _{n}$ as before to be the automorphism such that $%
\sigma _{n}(r_{n})=t_{n}$ and $\sigma _{n}^{2}=1$. Now, let $w_{n}$ be the
canonical unitary of $C^{\ast }\left( r_{n},t_{n}\right) \rtimes _{\sigma
_{n}}\mathbb{Z}_{2}$. Then we have the following relations:%
\begin{equation}
\left\{ 
\begin{array}{l}
w_{n}^{2}=1\text{,} \\ 
w_{n}t_{n}=r_{n}w_{n}\text{,} \\ 
t_{n}^{\ast }t_{n}=r_{n}^{\ast }r_{n}=t_{n}t_{n}^{\ast }+r_{n}r_{n}^{\ast }=1%
\end{array}%
\right.  \label{Cuntz-n}
\end{equation}%
Now, the fixed point C*-subalgebra of $C^{\ast }\left( r_{n},t_{n}\right) $
for $\sigma _{n}$ is generated by:%
\begin{equation}
\left\{ 
\begin{array}{l}
r_{n-1}=\frac{1}{\sqrt[2]{2}}\left( r_{n}+t_{n}\right) \text{,} \\ 
t_{n-1}=w_{n-1}\left( r_{n-1}\right) w_{n-1}\text{,}%
\end{array}%
\right.  \label{fixed}
\end{equation}%
where:%
\begin{equation}
w_{n-1}=r_{n}r_{n}^{\ast }-t_{n}t_{n}^{\ast }\text{.}  \label{unitaryformula}
\end{equation}%
By Proposition (\ref{decomposition1}), we have $w_{n-1}^{2}=1$ and $%
w_{n-1}t_{n-1}=r_{n-1}w_{n-1}$. By Theorem (\ref{FixedPointAlgebraThm}) we
have $r_{n-1}$ and $t_{n-1}$ are isometries, such that $r_{n-1}r_{n-1}^{\ast
}+t_{n-1}t_{n-1}^{\ast }=1$. Hence Relations (\ref{Cuntz-n}) are satisfied
for $n-1$. Therefore, $C^{\ast }\left( r_{n},t_{n}\right) $ is the
C*-crossed-product of $C^{\ast }\left( r_{n-1},t_{n-1}\right) $ for the
action of $\mathbb{Z}_{2}$ generated by $\sigma _{n-1}$ where $\sigma
_{n-1}(r_{n-1})=t_{n-1}$ and $\sigma _{n-1}^{2}=1$.

\bigskip Now, it is natural to define:%
\begin{equation*}
\left\{ 
\begin{array}{l}
r_{n+1}=\frac{1}{\sqrt[2]{2}}\left( r_{n}+r_{n}w_{n}\right) \text{,} \\ 
t_{n+1}=\frac{1}{\sqrt[2]{2}}\left( r_{n}-r_{n}w_{n}\right) \text{.}%
\end{array}%
\right.
\end{equation*}%
Thus, by Theorem (\ref{CrossedProductThm}), we have that $C^{\ast }\left(
r_{n},t_{n}\right) \rtimes _{\sigma _{n}}\mathbb{Z}_{2}=C^{\ast }\left(
r_{n},t_{n},w_{n}\right) $ is $C^{\ast }\left( r_{n+1},t_{n+1}\right) $.
Moreover, we note:%
\begin{equation*}
w_{n}=r_{n+1}r_{n+1}^{\ast }-t_{n+1}t_{n+1}
\end{equation*}%
which is of course Equation (\ref{unitaryformula}). Moreover, one checks
easily that Relation (\ref{fixed}) is satisfied for $n$ rather than $n-1$:%
\begin{equation*}
\frac{1}{\sqrt[2]{2}}\left( r_{n+1}+t_{n+1}\right) =\frac{1}{2}\left(
2r_{n}\right) =r_{n}\text{,}
\end{equation*}%
and%
\begin{eqnarray*}
\frac{1}{\sqrt[2]{2}}w_{n}\left( r_{n+1}+t_{n+1}\right) w_{n} &=&\frac{1}{2}%
\left[ w_{n}r_{n}w_{n}+w_{n}r_{n}+w_{n}r_{n}w_{n}-w_{n}r_{n}\right] \\
&=&\frac{1}{2}\left[ t_{n}+t_{n}\right] =t_{n}\text{ since }%
w_{n}r_{n}w_{n}=t_{n}\text{.}
\end{eqnarray*}

Thus, in particular, $C^{\ast }\left( t_{n},r_{n}\right) $ is the fixed
point C*-subalgebra of $C^{\ast }\left( t_{n+1},r_{n+1}\right) $ for the
action of $\mathbb{Z}_{2}$ generated by $\sigma _{n+1}$ which switches the
two generators $t_{n+1}$ and $r_{n+1}$. Thus, we have a pattern repeating
for $n\in \mathbb{Z}$ where $\mathcal{O}_{2}$ embeds in $\mathcal{O}_{2}$
either as a fixed point C*-subalgebra for the action which exchanges a
choice of generators of the target $\mathcal{O}_{2}$ or as the
crossed-product of the source $\mathcal{O}_{2}$ by the action which
exchanges a corresponding choice of generators of the source $\mathcal{O}%
_{2} $. Once a particular set of generators is chosen in our sequence, then
all the other ones are determined uniquely. Note that by Proposition (\ref%
{uniqueness}), the operators $w_{n}$ ($n\in \mathbb{Z}$) are then unique up
to a sign as well.

\section{Crossed-Product with $\mathbb{Z}$}

\bigskip Our study of the crossed-product $\mathcal{O}_{2}\rtimes _{\sigma }%
\mathbb{Z}_{2}$ makes it easy to study the crossed-product $\mathcal{O}%
_{2}\rtimes _{\sigma }\mathbb{Z}$. We now present a description of $\mathcal{%
O}_{2}\rtimes _{\sigma }\mathbb{Z}$. We begin with a simple observation:

\begin{proposition}
\label{IrredZ}Let $S_{1},S_{2}$ be two isometries with $S_{1}S_{1}^{\ast
}+S_{2}S_{2}^{\ast }=1$. Hence, $C^{\ast }\left( S_{1},S_{2}\right) =%
\mathcal{O}_{2}$. Let $\sigma $ be the automorphism defined by $\sigma
(S_{1})=S_{2}$ and $\sigma \left( S_{2}\right) =S_{1}$. Let $\pi $ be an
irreducible representation of $C^{\ast }\left( S_{1},S_{2}\right) \rtimes
_{\sigma }\mathbb{Z}$ and let $V$ be the canonical unitary in $C^{\ast
}\left( S_{1},S_{2}\right) \rtimes _{\sigma }\mathbb{Z}$. Then there exists $%
t\in \left[ -1,1\right] $ such that:%
\begin{equation*}
\pi (V)=e^{i\frac{\pi }{2}t}W
\end{equation*}%
with $W^{2}=1$ and $W\pi (S_{1})=\pi (S_{2})W$.
\end{proposition}

\begin{proof}
By construction, $V^{2}$ is in the center of $\mathcal{O}_{2}\rtimes
_{\sigma }\mathbb{Z}$. Since $\pi $ is irreducible, $\pi (V^{2})$ is a
scalar of the form $e^{i\pi t}$ for some $t\in \left[ -1,1\right] $. Let $%
W=e^{-i\frac{\pi }{2}t}V$. Then by construction, $W$ is a unitary such that $%
W^{2}=1$ and $W\pi (S_{1})W=V\pi (S_{1})V^{\ast }=\pi (S_{2})$ as desired.
\end{proof}

\bigskip We now can derive the following theorem:

\begin{theorem}
\label{ZCross}Let $S_{1},S_{2}$ be two isometries with $S_{1}S_{1}^{\ast
}+S_{2}S_{2}^{\ast }=1$. Hence, $C^{\ast }\left( S_{1},S_{2}\right) =%
\mathcal{O}_{2}$. Let $\sigma $ be the automorphism defined by $\sigma
(S_{1})=S_{2}$ and $\sigma \left( S_{2}\right) =S_{1}$. Then $\mathcal{O}%
_{2}\rtimes _{\sigma }\mathbb{Z}$ is *-isomorphic to:%
\begin{equation*}
\left\{ f\in C\left( \left[ -1,1\right] ,\mathcal{O}_{2}\rtimes _{\sigma }%
\mathbb{Z}_{2}\right) :f(-1)=\widehat{\sigma }(f(1))\right\} .
\end{equation*}
\end{theorem}

\begin{proof}
To fix notation, let us write:%
\begin{equation*}
\mathcal{A}=\left\{ f\in C\left( \left[ -1,1\right] ,\mathcal{O}_{2}\rtimes
_{\sigma }\mathbb{Z}_{2}\right) :f(-1)=\widehat{\sigma }(f(1))\right\}
\end{equation*}%
where we denote by $w$ the canonical unitary in $\mathcal{O}_{2}\rtimes
_{\sigma }\mathbb{Z}_{2}$ and $\widehat{\sigma }$ is the automorphism
defined uniquely by $\widehat{\sigma }(a)=a$ for $a\in \mathcal{O}_{2}$ and $%
\widehat{\sigma }(w)=-w$.

We introduce the following elements of $\mathcal{A}$:%
\begin{equation*}
\left\{ 
\begin{array}{l}
v:t\in \left[ -1,1\right] \mapsto e^{i\pi \frac{t}{2}}w\text{,} \\ 
s_{1}:t\in \left[ -1,1\right] \mapsto S_{1}\text{,} \\ 
s_{2}:t\in \left[ -1,1\right] \mapsto S_{2}\text{.}%
\end{array}%
\right.
\end{equation*}

Our proof consists of two steps:\ we show first that $\mathcal{A}=C^{\ast
}\left( s_{1},s_{2},v\right) $. We then show that $C^{\ast }\left(
s_{1},s_{2},v\right) $ is *-isomorphic to $\mathcal{O}_{2}\rtimes _{\sigma }%
\mathbb{Z}$.

By construction, $C^{\ast }\left( s_{1},s_{2},v\right) \subseteq \mathcal{A}$%
. To show that $\mathcal{A}=C^{\ast }\left( s_{1},s_{2},v\right) $, we
introduce the C*-subalgebra $\mathcal{B}$ of $\mathcal{A}$ defined by:%
\begin{equation*}
\left\{ f\in C\left( \left[ -1,1\right] ,\mathcal{O}_{2}\rtimes _{\sigma }%
\mathbb{Z}^{2}\right) :f(-1)=f(1)\right\} \text{.}
\end{equation*}%
Writing $f=\frac{1}{2}\left( f+\widehat{\sigma }(f)\right) +\frac{1}{2}%
\left( f-\widehat{\sigma }(f)\right) v^{2}$ (since $v^{2}=1$), we easily see
that:%
\begin{equation*}
\mathcal{A}=\left\{ f+gv:f,g\in \mathcal{B}\right\} \text{.}
\end{equation*}%
Thus, to prove $\mathcal{A}\subseteq C^{\ast }\left( s_{1},s_{2},v\right) $
it is enough to show that$\mathcal{\ B}\subseteq C^{\ast }\left(
s_{1},s_{2},v\right) $. Now, $s_{1},s_{2}$ and $v^{2}:t\in \left[ -1,1\right]
\mapsto e^{i\pi t}$ are all in $\mathcal{B}$ by construction, and a standard
argument shows that:%
\begin{equation*}
\mathcal{B}=C^{\ast }\left( s_{1},s_{2},v^{2}\right) \cong \mathcal{O}%
_{2}\otimes C\left( \left[ -1,1\right] \right)
\end{equation*}%
so $\mathcal{B}\subseteq C^{\ast }\left( s_{1},s_{2},v\right) $.

Now, it is sufficient to show that $C^{\ast }\left( s_{1},s_{2},v\right) $
is *-isomorphic to $\mathcal{O}_{2}\rtimes _{\sigma }\mathbb{Z}$. Let $V$ be
the canonical unitary of $\mathcal{O}_{2}\rtimes _{\sigma }\mathbb{Z}$. \
Since $vs_{1}v^{\ast }=s_{2}$ and $vs_{2}v^{\ast }=s_{1}$ by construction,
there exists by universality of the crossed-product a unique *-epimorphism $%
\theta :\mathcal{O}_{2}\rtimes _{\sigma }\mathbb{Z}\longrightarrow C^{\ast
}\left( s_{1},s_{2},v\right) $ with $\theta (S_{1})=s_{1}$, $\theta
(S_{2})=s_{2}$ and $\theta (V)=v$. We wish to show that $\theta $ is in fact
a *-isomorphism. Let $a\in \ker \theta $. Let $\pi $ be an arbitrary
irreducible representation of $\mathcal{O}_{2}\rtimes _{\sigma }\mathbb{Z}$.
By Proposition (\ref{IrredZ}), there exists $t\in \left[ -1,1\right] $ such
that $\pi (V)=e^{i\frac{\pi }{2}t}W$ with $W^{2}=1$ and $W\pi (S_{1})=\pi
(S_{2})W$ with $W$ unitary. By universality, there exists a representation $%
\psi $ of $\mathcal{O}_{2}\rtimes _{\sigma }\mathbb{Z}_{2}$ on the same
Hilbert space on which $\pi $ acts such that $\psi (S_{1})=\pi (S_{1}),$ $%
\psi (S_{2})=\pi (S_{2})$ and $\psi (w)=W$. Let $\varepsilon _{t}$ be the
*-morphism from $\mathcal{A}$ onto $\mathcal{O}_{2}\rtimes _{\sigma }\mathbb{%
Z}_{2}$ defined by $\varepsilon _{t}(f)=f(t)$ for all $f\in \mathcal{A}$.
Then by construction, $\pi =\psi \circ \varepsilon _{t}\circ \theta $. Hence 
$\pi (a)=0$. Since $\pi $ was arbitrary irreducible, $a=0$ and thus $\theta $
is injective. This completes our proof.
\end{proof}

\begin{corollary}
Let $S_{1},S_{2}$ be two isometries with $S_{1}S_{1}^{\ast
}+S_{2}S_{2}^{\ast }=1$. Hence, $C^{\ast }\left( S_{1},S_{2}\right) =%
\mathcal{O}_{2}$. Let $\sigma $ be the automorphism defined by $\sigma
(S_{1})=S_{2}$ and $\sigma \left( S_{2}\right) =S_{1}$. Then $\mathcal{O}%
_{2}\rtimes _{\sigma }\mathbb{Z}$ is *-isomorphic to:%
\begin{equation*}
\left\{ f\in C\left( \left[ -1,1\right] ,\mathcal{O}_{2}\right)
:f(-1)=\sigma (f(1))\right\} \text{.}
\end{equation*}
\end{corollary}

\begin{proof}
By Proposition (\ref{Tau}), there exists a *-isomorphism:%
\begin{equation*}
\tau :\mathcal{O}_{2}\rtimes \mathbb{Z}_{2}\longrightarrow \mathcal{O}_{2}
\end{equation*}%
such that $\sigma \circ \tau =\tau \circ \widehat{\sigma }$, where we use
the notations in the proof of Theorem (\ref{ZCross}).

The corollary follows from this observation and Theorem (\ref{ZCross}).
\end{proof}

We can rephrase the result above in a manner which may appear explicit. We
call an element $a$ of $\mathcal{O}_{2}$ \emph{symmetric} if $a=\sigma (a)$
and \emph{antisymmetric} if $a=-\sigma (a)$. Then we get immediately from
Theorem (\ref{ZCross}):

\begin{corollary}
We have:%
\begin{equation*}
\mathcal{O}_{2}\rtimes _{\sigma }\mathbb{Z}=\left\{ f\in C\left( \left[ -1,1%
\right] ,\mathcal{O}_{2}\right) ~~\left\vert 
\begin{array}{l}
f(1)+f(-1)\text{ is symmetric,} \\ 
f(1)-f(-1)\text{ is antisymmetric.}%
\end{array}%
\right. \right\}
\end{equation*}
\end{corollary}

\section{Appendix:\ Concrete Irreducible Representations}

\bigskip In this appendix, we present a concrete representation of $\mathcal{%
O}_{2}$ which fits the framework of this paper. Our representation is based
upon the following group:

\begin{definition}
Let $\mathcal{A}$ be the group of strictly increasing affine transformations
of $\mathbb{R}$, i.e.:%
\begin{equation*}
\mathcal{A=}\left\{ \varphi _{a,b}:t\in \mathbb{R}\mapsto at+b:a>0,b\in 
\mathbb{R}\right\} \text{.}
\end{equation*}
\end{definition}

The group $\mathcal{A}$ is naturally isomorphic to $\left\{ \left[ 
\begin{array}{cc}
1 & b \\ 
0 & a%
\end{array}%
\right] :a>0,b\in \mathbb{R}\right\} $ where $\left[ 
\begin{array}{cc}
1 & b \\ 
0 & a%
\end{array}%
\right] $ is mapped to $\varphi _{a,b}$. We will use this isomorphism
implicitly when convenient.

\begin{proposition}
For any $\varphi _{a,b}\in \mathcal{A}$ we define the bounded linear
operator $\pi _{\varphi _{a,b}}$ of $L^{2}\left( \mathbb{R}\right) $ by:%
\begin{equation*}
\pi _{\varphi _{a,b}}:f\in L^{2}\left( \mathbb{R}\right) \mapsto a^{\frac{1}{%
2}}f\circ \varphi _{a,b}\text{.}
\end{equation*}%
Then $\pi $ is a unitary representation of $\mathcal{A}$ on $L^{2}\left( 
\mathbb{R}\right) $.
\end{proposition}

\begin{proof}
It is immediate that $f\mapsto f\circ \varphi _{a,b}$ is a linear operator
on $L^{2}\left( \mathbb{R}\right) $ and $\pi _{gg^{\prime }}=\pi _{g}\pi
_{g^{\prime }}$ for all $g,g^{\prime }\in \mathcal{A}$. Moreover $\pi _{%
\limfunc{id}}=\limfunc{id}$. Now for all $f,g\in L^{2}\left( \mathbb{R}%
\right) $ we have:%
\begin{equation*}
\int_{\mathbb{R}}f\left( at+b\right) g(t)dt=\int_{\mathbb{R}}f\left(
t\right) \frac{1}{a}g\left( \frac{1}{a}\left( t-b\right) \right) dt
\end{equation*}%
so $\left\langle \pi _{\varphi _{a,b}}\left( f\right) ,g\right\rangle
=\left\langle f,\pi _{\varphi _{\frac{1}{a},-\frac{b}{a}}}\left( g\right)
\right\rangle =\left\langle f,\pi _{\varphi _{a,b}^{-1}}(g)\right\rangle $
and thus $\pi _{\varphi _{a,b}}$ is bounded and unitary.
\end{proof}

\bigskip

\begin{definition}
Let $I$ be any closed subset in $\mathbb{R}$. The orthogonal projection from 
$L^{2}\left( \mathbb{R}\right) $ onto $L^{2}\left( I\right) $ is denoted by $%
P_{I}$.
\end{definition}

In other words, $P_{I}$ is the multiplication operator by the indicator
function $\chi _{I}$ of $I$.

\begin{definition}
\label{A1D1}Let $I,J$ be two compact intervals in $\mathbb{R}$. Let $\varphi\in \mathcal{A}$ be the unique increasing affine map such that $%
\varphi (I)=J$. Then we set:%
\begin{equation*}
V\left( I,J\right) =P_{I}\pi _{\varphi }P_{J}\text{.}
\end{equation*}
\end{definition}

Note that $P_{I}\pi =P_{I}\pi _{\varphi }P_{J}=\pi _{\varphi }P_{J}$ by
construction in Definition (\ref{A1D1}).

\begin{theorem}
\label{Appendix1}The set:%
\begin{equation*}
\Sigma =\left\{ 0\right\} \cup \left\{ V\left( I,J\right) :I,J\text{ compact
intervals in }\mathbb{R}\right\}
\end{equation*}%
is a semigroup of partial isometries. Moreover:

\begin{enumerate}
\item For all compact interval $I$ we have $V(I,I)=P_{I}$,

\item For all compact intervals $I,J$ we have $V(I,J)=V(J,I)^{\ast }$,

\item For all four compact intervals $I,J,K,L$ we have:%
\begin{equation*}
V(I,J)V(K,L)=V\left( \varphi _{1}^{-1}\left( J\cap K\right) ,\varphi
_{2}\left( J\cap K\right) \right) \text{.}
\end{equation*}%
where $\varphi _{1}$ and $\varphi _{2}$ are the unique elements of $\mathcal{%
A}$ such that $\varphi _{1}(I)=J$ and $\varphi _{2}(K)=L$. Note that in
particular:%
\begin{equation*}
\varphi _{\varphi _{1}^{-1}\left( J\cap K\right) ,\varphi _{2}\left( J\cap
K\right) }=\varphi _{2}\varphi _{1}\text{.}
\end{equation*}
\end{enumerate}

In particular, the initial space of $V(I,J)$ is $L^{2}(J)$ and the final
space is $L^{2}(I)$ for all compact intervals $I,J$ of $\mathbb{R}$.
\end{theorem}

\begin{proof}
By uniqueness of the element in $\mathcal{A}$ which maps an interval to
another, properties 1 and 2 are immediate. In general, given four compact
intervals $J_{1},J_{2},J_{3}$ and $J_{4}$, and two affine maps $\varphi _{1}$
and $\varphi _{2}$ such that $\varphi \left( J_{1}\right) =J_{2}$ and $%
\varphi ^{\prime }\left( J_{3}\right) =J_{4}$, then we let $J=J_{2}\cap
J_{3} $.%
\begin{eqnarray*}
V\left( J_{1},J_{2}\right) V\left( J_{3},J_{4}\right) &=&P_{J_{1}}\pi
_{\varphi _{1}}P_{J_{2}}P_{J_{3}}\pi _{\varphi _{2}}P_{J_{4}} \\
&=&P_{J_{1}}\pi _{\varphi _{1}}P_{J_{2}\cap J_{3}}\pi _{\varphi
_{2}}P_{J_{4}}\text{.}
\end{eqnarray*}%
Now, for $\pi _{\varphi _{2}}\left( f\right) $ to be supported on $J_{2}\cap
J_{3}$ we must have that $f$ is supported on $\varphi _{2}\left( J_{2}\cap
J_{3}\right) $. Hence:%
\begin{equation*}
P_{J_{2}\cap J_{3}}\pi _{\varphi _{2}}(f)P_{J_{4}}=P_{J_{2}\cap J_{3}}\pi
_{\varphi _{2}}(f)P_{\varphi _{2}\left( J_{2}\cap J_{3}\right) }\text{.}
\end{equation*}

Similarly, if $f$ is supported on $J_{2}\cap J_{3}$ then $\pi _{\varphi
_{1}}(f)$ is supported on $\varphi _{1}^{-1}\left( J_{2}\cap J_{3}\right) $
and:%
\begin{equation*}
P_{J_{1}}\pi _{\varphi _{1}}P_{J_{2}\cap J_{3}}=P_{\varphi _{1}^{-1}\left(
J_{2}\cap J_{3}\right) }\pi _{\varphi _{1}}P_{J_{2}\cap J_{3}}\text{.}
\end{equation*}%
Hence the third property above.
\end{proof}

\bigskip We will use two simple lemmas to prove that the representation of $
\mathcal{O}_{2}$ introduced in Theorem (\ref{O2representation}) is irreducible. Let $\lambda $ be the usual Lebesgue measure on $\left[ -1,1\right] $ (so that $\lambda \left( \left[ -1,1\right]
\right) =2$).

\begin{lemma}
\label{AppLemma1}Let $\mathfrak{a}$ be an arbitrary Borel subset of $\left(
-1,1\right) $ of strictly positive Lebesgue measure in $\left( 0,2\right) $. Then there exist a natural number $m$ and two integers $k_{1}$ and $k_{2}$ such that 
\begin{equation*}
\lambda \left( \mathfrak{a}\cap \left[ \frac{k_{1}}{2^{m}},\frac{k_{1}+1}{%
2^{m}}\right] \right) >\lambda \left( \mathfrak{a}\cap \left[ \frac{k_{2}}{%
2^{m}},\frac{k_{2}+1}{2^{m}}\right] \right) 
\end{equation*}%
with $-2^{m}\leq k_{1},k_{2}<2^{m}$.
\end{lemma}

\begin{proof}
Since $\lambda \left( \mathfrak{a}\right) <2$, there exists an open set $%
\mathfrak{g}$ in $\left( -1,1\right) $ such that $\mathfrak{a}\subseteq 
\mathfrak{g}$ and $\lambda \left( \mathfrak{g}\right) <2$. Now, $\mathfrak{g}
$ is the disjoint union countably many open intervals $\mathfrak{g}_{i}$ $%
\left( i\in \mathbb{N}\right) $ in $\left[ -1,1\right] $. In fact we may
choose each $\mathfrak{g}_{i}$ of the form $\left( \frac{k_{i}}{2^{m_{i}}},%
\frac{k_{i}+1}{2^{m_{i}}}\right) $ for some $k_{i}\in \mathbb{Z}$ and $%
m_{i}\in \mathbb{N}$ for all $i\in \mathbb{N}$ -- in which case the
symmetric difference between $\mathfrak{g}$ and $\bigcup_{i\in 
\mathbb{N}}\mathfrak{g}_{i}$ has measure $0$. Now:%
\begin{equation*}
\sum_{i\in \mathbb{N}}\lambda \left( \mathfrak{g}_{i}\cap \mathfrak{a}%
\right) =\lambda \left( \mathfrak{a}\right) <\lambda \left( \mathfrak{a}%
\right) \frac{1}{2}\lambda \left( \mathfrak{g}\right) =\frac{\lambda \left( 
\mathfrak{a}\right) }{2}\sum_{i\in \mathbb{N}}\lambda \left( \mathfrak{g}%
_{i}\right) 
\end{equation*}%
so there exists $i\in \mathbb{N}$ such that $\lambda \left( \mathfrak{g}%
_{i}\cap \mathfrak{a}\right) >\frac{\lambda \left( \mathfrak{a}\right) }{2}%
\lambda \left( \mathfrak{g}_{i}\right) $. To fix notations, let us write $%
\mathfrak{g}_{i}=\left[ \frac{k_{1}}{2^{m}},\frac{k_{1}+1}{2^{m}}\right] $,
so that:%
\begin{equation*}
\lambda \left( \mathfrak{a}\cap \left[ \frac{k_{1}}{2^{m}},\frac{k_{1}+1}{%
2^{m}}\right] \right) >\frac{1}{2^{m+1}}\lambda \left( \mathfrak{a}\right) 
\text{.}
\end{equation*}

On the other hand:%
\begin{equation*}
\sum_{k=-2^{m}}^{2^{m}-1}\lambda \left( \mathfrak{a}\cap \left[ \frac{k}{%
2^{m}},\frac{k+1}{2^{m}}\right] \right) =\lambda \left( \mathfrak{a}\right)
\end{equation*}%
so there exists an integer $k_{2}$ such that:%
\begin{equation*}
\lambda \left( \mathfrak{a}\cap \left[ \frac{k_{2}}{2^{m}},\frac{k_{2}+1}{%
2^{m}}\right] \right) \leq \frac{1}{2^{m+1}}\lambda \left( \mathfrak{a}%
\right) \text{.}
\end{equation*}%
Consequently:%
\begin{equation*}
\lambda \left( \mathfrak{a}\cap \left[ \frac{k_{2}}{2^{m}},\frac{k_{2}+1}{%
2^{m}}\right] \right) <\lambda \left( \mathfrak{a}\cap \left[ \frac{k_{1}}{%
2^{m}},\frac{k_{1}+1}{2^{m}}\right] \right)
\end{equation*}%
as desired.
\end{proof}

\begin{lemma}
\label{AppLemma2}Let $J_{1}$ and $J_{2}$ be two closed intervals in $\left[
-1,1\right] $. Let $\mathfrak{a}$ be a Borel subset of $\left[ -1,1\right] $
such that the projection $P_{\mathfrak{a}}$ commutes with $V\left(
J_{2},J_{1}\right) $. Then $\lambda \left( J_{1}\cap \mathfrak{a}\right)
=\lambda \left( J_{2}\cap \mathfrak{a}\right) $.
\end{lemma}

\begin{proof}
Write $J=J_{1}$ and define $c\in \left[ -1,1\right] $ by $J_{2}=J_{1}+c$.
Write $P=P_{\mathfrak{a}}$ and $V=V\left( J+c,J\right) $. Now for $f\in
L^{2}\left( \left[ -1,1\right] \right) $ and $t\in \left[ -1,1\right] $ we
have:%
\begin{equation*}
Pf(t)=\chi _{\mathfrak{a}}(t)f(t)\text{ and }Vf\left( t\right) =\chi
_{J+c}(t)f(t-c)\text{.}
\end{equation*}%
Thus $PV=VP$ exactly when:%
\begin{equation*}
\chi _{\mathfrak{a}}(t)\chi _{J+c}(t)f(t-c)=\chi _{J+c}(t)\chi _{\mathfrak{a}%
}(t-c)f(t-c)
\end{equation*}%
for all $f\in L^{2}\left( \left[ -1,1\right] \right) $ and $t\in \left[ -1,1%
\right] $. Thus $\chi _{\mathfrak{a}\cap \left( J+c\right) }=\chi _{\left( 
\mathfrak{a}\cap J\right) +c}$, which implies the desired result.
\end{proof}

\begin{theorem}\label{O2representation}
Let:%
\begin{eqnarray*}
&&S_{1}\text{ be the restriction of }V\left( \left[ 0,1\right] ,\left[ -1,1%
\right] \right) \text{ to }L^{2}\left( \left[ -1,1\right] \right) \text{,} \\
&&S_{2}\text{ be the restriction of }V\left( \left[ -1,0\right] ,\left[ -1,1%
\right] \right) \text{ to }L^{2}\left( \left[ -1,1\right] \right) \text{.}
\end{eqnarray*}%
In other words, for $f\in L^{2}\left( \left[ -1,1\right] \right) $ and $t\in %
\left[ -1,1\right] $ we have:%
\begin{equation*}
S_{1}\left( f\right) \left( t\right) =\sqrt[2]{2}f\left( 2t-1\right) \text{
and }S_{2}\left( f\right) \left( t\right) =\sqrt[2]{2}f\left( 2t+1\right) 
\text{.}
\end{equation*}

Then $S_{1}$ and $S_{2}$ are two isometries of $L^{2}\left( \left[ -1,1%
\right] \right) $ such that $S_{1}S_{1}^{\ast }+S_{2}S_{2}^{\ast }=1$.
Moreover $C^{\ast }\left( S_{1},S_{2}\right) $ is irreducible and:%
\begin{equation*}
C^{\ast }(S_{1},S_{2})=\overline{\limfunc{span}\left\{ V\left( I,J\right)
:I,J\in \mathcal{J}\right\} }
\end{equation*}%
where:%
\begin{equation*}
\mathcal{J}=\left\{ \left[ \frac{k}{2^{m}},\frac{k^{\prime }}{2^{m}}\right]
:m\in \mathbb{N}\text{, }-2^{m}\leq k<k^{\prime }\leq 2^{m}\text{, }%
k,k^{\prime }\text{ integers}\right\} \text{.}
\end{equation*}
\end{theorem}

\begin{proof}
By construction, $S_{1}$ and $S_{2}$ are isometries. Moreover:%
\begin{equation*}
S_{1}S_{1}^{\ast }=P_{\left[ 0,1\right] }\text{ and }S_{2}S_{2}^{\ast }=P_{%
\left[ -1,0\right] }
\end{equation*}%
so $S_{1}S_{1}^{\ast }+S_{2}S_{2}^{\ast }=1$ in $L^{2}\left( \left[ -1,1%
\right] \right) $.

Note that by definition, the semigroup generated by $S_{1}$ and $S_{2}$ is
the semigroup generated by $V\left( \left[ 0,1\right] ,\left[ -1,1\right]
\right) $ and $V\left( \left[ -1,0\right] ,\left[ -1,1\right] \right) $ when
regarded as operators acting on $L^{2}\left( \left[ -1,1\right] \right) $
only. First, we observe that $\varphi _{S_{1}}:t\mapsto 2t+1$ and $\varphi
_{S_{2}}:t\mapsto 2t-1$. Hence, $\varphi _{S_{1}}^{-1}:t\mapsto \frac{1}{2}t-%
\frac{1}{2}$ and $\varphi _{S_{2}}^{-1}:t\mapsto \frac{1}{2}t+\frac{1}{2}$
and thus, by Theorem (\ref{Appendix1}):%
\begin{equation*}
V\left( \left[ 0,1\right] ,\left[ -1,1\right] \right) V\left( \left[ \frac{k%
}{2^{m}},\frac{k^{\prime }}{2^{m}}\right] ,\left[ -1,1\right] \right)
=V\left( \left[ \frac{k-1}{2^{m+1}},\frac{k^{\prime }-1}{2^{m+1}}\right] ,%
\left[ -1,1\right] \right) \text{.}
\end{equation*}%
and:%
\begin{equation*}
V\left( \left[ -1,0\right] ,\left[ -1,1\right] \right) V\left( \left[ \frac{k%
}{2^{m}},\frac{k^{\prime }}{2^{m}}\right] ,\left[ -1,1\right] \right)
=V\left( \left[ \frac{k+1}{2^{m+1}},\frac{k^{\prime }+1}{2^{m+1}}\right] ,%
\left[ -1,1\right] \right) \text{.}
\end{equation*}%
Thus by induction, any finite product of $S_{1}$ and $S_{2}$ is of the form $%
V\left( I,\left[ -1,1\right] \right) $ where $I\in \mathcal{J}$ and moreover
all such operators can be obtained as such finite products. So the semigroup
generated \ by $S_{1}$ and $S_{2}$ is given by:%
\begin{equation*}
\mathcal{S=}\left\{ V\left( I,\left[ -1,1\right] \right) :I\in \mathcal{J}%
\right\} \text{.}
\end{equation*}%
Now, by Theorem (\ref{Appendix1}), we also have that the adjoint of the
operators in $\mathcal{S}$ are of the form:%
\begin{equation*}
\mathcal{S}^{\ast }\mathcal{=}\left\{ V\left( \left[ -1,1\right] ,I\right)
:I\in \mathcal{J}\right\} \text{.}
\end{equation*}%
Hence by a direct computation and applying Theorem (\ref{Appendix1}), we get
that arbitrary products of $S_{1}$,$S_{1}^{\ast }$,$S_{2}$,$S_{2}^{\ast }$
are exactly given by $0$ or:%
\begin{equation*}
V\left( \left[ \frac{k}{2^{m}},\frac{k^{\prime }}{2^{m}}\right] ,\left[ 
\frac{p}{2^{m}},\frac{q}{2^{m}}\right] \right) =V\left( \left[ \frac{k}{2^{m}%
},\frac{k^{\prime }}{2^{m}}\right] ,\left[ -1,1\right] \right) V\left( \left[
-1,1\right] ,\left[ \frac{p}{2^{m}},\frac{q}{2^{m}}\right] \right)
\end{equation*}%
for all:%
\begin{equation*}
-2^{m}\leq p,k<q,k^{\prime }\leq 2^{m}\text{, }p,q,k,k^{\prime }\in \mathbb{Z%
}\text{ and }m\in \mathbb{N}\text{.}
\end{equation*}%
Therefore:%
\begin{equation*}
C^{\ast }(S_{1},S_{2})=\overline{\limfunc{span}\left\{ V\left( I,J\right)
:I,J\in \mathcal{J}\right\} }
\end{equation*}%
as claimed.

Last, note that $C^{\ast }\left( S_{1},S_{2}\right) $ contains $P_{I}$ for
all compact intervals $I$ with dyadic end points. Hence, the commutant of $%
C^{\ast }\left( S_{1},S_{2}\right) $ is contained in the Von Neumann algebra
of the multiplications operators by functions in $L^{\infty }\left( \left[
-1,1\right] \right) $ on $L^{2}\left( \left[ -1,1\right] \right) $.
Consequently, if $P$ is a projection commuting with $C^{\ast }\left(
S_{1},S_{2}\right) $ then there exists a measurable set $\mathfrak{a}%
\subseteq \left[ -1,1\right] $ such that $P$ is the multiplication operator
with the indicator function $\chi _{\mathfrak{a}}$ of $A$. Let us assume
that $\lambda (A)\in \left( 0,2\right) $. Then by Lemma (\ref{AppLemma1})\
we can find a natural number $m$ and two integers $k_{1}$ and $k_{2}$ such
that:%
\begin{equation*}
\lambda \left( \mathfrak{a}\cap \left[ \frac{k_{1}}{2^{m}},\frac{k_{1}+1}{%
2^{m}}\right] \right) >\lambda \left( \mathfrak{a}\cap \left[ \frac{k_{2}}{%
2^{m}},\frac{k_{2}+1}{2^{m}}\right] \right) \text{.}
\end{equation*}

Yet, this contradicts Lemma (\ref{AppLemma2}). So $\lambda \left( \mathfrak{a%
}\right) \in \left\{ 0,2\right\} $ and thus our representation is
irreducible.
\end{proof}

\bigskip We now can construct a unitary which implements the action of $%
\mathbb{Z}_{2}$ which flips $S_{1}$ and $S_{2}$ and illustrate our work by
applying Theorems\ (\ref{CrossedProductThm}) and (\ref{FixedPointAlgebraThm}%
) to describe the fixed point subalgebra and the crossed-product in term of
concrete operators.

\begin{remark}
Let $W:L^{2}\left( \left[ -1,1\right] \right) \rightarrow L^{2}\left( \left[
-1,1\right] \right) $ be defined by $Wf:t\in \left[ -1,1\right] \mapsto
f(-t) $ for all $f\in L^{2}\left( \left[ -1,1\right] \right) $. Then $W$ is
an order two unitary such that $WS_{2}=S_{1}W$. Thus, as in Theorem (\ref%
{CrossedProductThm}), setting $\sigma (A)=WAW$ for $A\in C^{\ast }\left(
S_{1},S_{2}\right) $, we conclude that $C^{\ast }\left( S_{1},S_{2}\right)
\rtimes _{\sigma }\mathbb{Z}_{2}=C^{\ast }\left( S_{1},S_{1}W\right) $ is
*-isomorphic to $\mathcal{O}_{2}$. Note that for $f\in L^{2}\left( \left[
-1,1\right] \right) $ and $t\in \left[ -1,1\right] $, and our choice of
representation in this section, we get $S_{1}(f)(t)=\sqrt[2]{2}f\left(
2t-1\right) $ and $\left( WS_{1}\right) \left( f\right) \left( t\right) =%
\sqrt[2]{2}f\left( -2t-1\right) $ -- thus both are isometries. One checks
easily that $S_{1}S_{1}^{\ast }+\left( WS_{1}\right) \left( WS_{1}\right)
^{\ast }=1$.
\end{remark}

\begin{remark}
We can describe the two generators of the fixed point C*-algebra of $C^{\ast
}\left( S_{1},S_{2}\right) $ given by Theorem (\ref{FixedPointAlgebraThm})\
in our concrete representation. Keeping the notations of Theorem (\ref%
{FixedPointAlgebraThm}), we have, for $f\in L^{2}\left( \left[ -1,1\right]
\right) $ and $t\in \left[ -1,1\right] $:%
\begin{eqnarray*}
U &=&S_{1}S_{1}^{\ast }-S_{2}S_{s}^{\ast }=P_{\left[ 0,1\right] }-P_{\left[
-1,0\right] }\text{,} \\
T &=&\sqrt[2]{2}\left( S_{1}+S_{2}\right) \text{ so }T\left( f\right) \left(
t\right) =f\left( 2t-1\right) +f\left( 2t+1\right) \\
R &=&UTU\text{ so }R\left( f\right) \left( t\right) =\left\{ 
\begin{array}{c}
f\left( 2t-1\right) \text{ for }t\in \left( \frac{1}{2},1\right) \text{,} \\ 
-f\left( 2t-1\right) \text{ for }t\in \left( 0,\frac{1}{2}\right) \text{,}
\\ 
-f\left( 2t+1\right) \text{ for }t\in \left( -\frac{1}{2},0\right) \text{,}
\\ 
f\left( 2t+1\right) \text{ for }t\in \left( -1,-\frac{1}{2}\right) \text{.}%
\end{array}%
\right.
\end{eqnarray*}%
One then check that $R$ and $T$ are isometries such that $TT^{\ast
}+RR^{\ast }=1$.
\end{remark}

\bibliographystyle{amsplain}
\bibliography{../thesis.bib}

\providecommand{\bysame}{\leavevmode\hbox to3em{\hrulefill}\thinspace}
\providecommand{\MR}{\relax\ifhmode\unskip\space\fi MR }
\providecommand{\MRhref}[2]{%
  \href{http://www.ams.org/mathscinet-getitem?mr=#1}{#2}
}
\providecommand{\href}[2]{#2}
\begin{thebibliography}{1}

\bibitem{Choi79}
{M}.-{D}. {C}hoi, \emph{A simple {$C^\ast$}-algebra generated by two
  finite-order unitaries}, Canad. J. Math. \textbf{31} (1979), no.~4, 867--880.

\bibitem{Latremoliere06b}
{M}.-{D}. {C}hoi and {F}. {L}atr{\'e}moli{\`e}re, \emph{The {C*}-algebra of
  symmetric words in two universal unitaries}, {J}. {O}perator {T}heory
  \textbf{62} (2009), 159--169, math.OA/0610467.

\bibitem{Cuntz77}
{J}. {C}untz, \emph{{S}imple {$C^\ast$}-algebras generated by isometries},
  Comm. Math. Phys. \textbf{57} (1977), 173--185.

\bibitem{Davidson}
{K}.~{R}. {D}avidson, \emph{{C*}--algebras by example}, Fields Institute
  Monographs, American Mathematical Society, 1996.

\bibitem{Izumi04}
{M}. {I}zumi, \emph{Finite group actions on {$C\sp\ast$}-algebras with the
  rohlin property {I}}, Duke Math. J. \textbf{122} (2004), no.~2, 233--280.

\bibitem{Izumi04b}
\bysame, \emph{Finite group actions on {$C\sp\ast$}-algebras with the rohlin
  property {II}}, Advances in Math. \textbf{184} (2004), 119--160.

\bibitem{Rieffel80}
{M}. {R}ieffel, \emph{Actions of finite groups on {C*}-algebras}, Math. Scand.
  \textbf{47} (1980), 157--176.

\bibitem{Rosenberg79}
{J}. {R}osenberg, \emph{Appendix to {O}. {B}ratteli paper on
  {``Crossed-products of UHF algebras"}}, Duke Math. J. \textbf{46} (1979),
  no.~1, 25--26.

\bibitem{Zeller-Meier68}
{G}. {Z}eller{-}{M}eier, \emph{Produits crois{\'e}s {d'u}ne {C*-}alg{\`e}bre
  par un groupe {d' A}utomorphismes}, {J}. {M}ath. pures et appl. \textbf{47}
  (1968), no.~2, 101--239.

\end{thebibliography}

\end{document}